\documentclass[a4paper, 10pt, conference]{ieeeconf}      

\IEEEoverridecommandlockouts                              
\usepackage{graphics,algorithm} 
\usepackage{epsfig} 
\usepackage{mathptmx} 
\usepackage{times} 
\usepackage{amsmath} 
\usepackage{amssymb}  

\usepackage[mathcal]{euscript} 
\usepackage{mathtools}
\usepackage{mathrsfs}

\usepackage{graphicx}
\graphicspath{{Figures/}} 
\usepackage{caption} 
\usepackage{subfig} 

\usepackage{url}

\usepackage{xcolor}
\definecolor{USred}{rgb}{0.74,0.1,0.1}
\definecolor{USblue}{rgb}{0.2,0.2,0.7}
\definecolor{green1}{cmyk}{0.82,0,1,0.3}

\definecolor{Royalblue}{cmyk}{1,0.30,0.2,0.2}

\newcommand{\al}[1]{\begin{align} #1\end{align}}
\newcommand{\nn}{\nonumber}
\newcommand{\argmin}{\operatornamewithlimits{argmin}}

\DeclareRobustCommand{\vect}[1]{
	\ifcat#1\relax
	\boldsymbol{#1}
	\else
	\mathbf{#1}
	\fi}

\newcommand{\col}[1]{\mathrm{col}(#1)}

\title{\LARGE \bf Distributed Kalman filtering with event-triggered communication:\\ a robust approach}

\author{Davide Ghion, Mattia Zorzi 
\thanks{D. Ghion is with Serenissima Informatica Spa, (e-mail: davide.ghion@serinf.it)}
\thanks{M. Zorzi is with the Department of Information Engineering, University of Padova, Padova, Italy; email:	 
	   {\tt\small zorzimat@dei.unipd.it}}%
\thanks{}%
}

\begin{document}

\maketitle
\thispagestyle{empty}
\pagestyle{empty}

\begin{abstract}       
We consider the problem of distributed Kalman filtering for sensor networks in the case there is a limit in data transmission  and there is model uncertainty. More precisely, we propose a distributed filtering strategy with event-triggered communication in which the state estimators are computed according to the least favorable model. The latter belongs to a ball (in Kullback-Leibler topology) about the nominal model. We also present a preliminary numerical example in order to test the performance of the proposed strategy.         
\end{abstract}

\section{Introduction}
Nowadays, sensor networks play an important role in various fields such as  security, monitoring, data analysis and so on.
In these applications, the sensors collect the measurements and from them it is required to estimate the state of a dynamic system. This task is performed in a distributed fashion and it can be accomplished in different ways, see \cite{cattivelli2010,spanos2005,Saber_CDC2007,BATTISTELLI2014707,6576197}; for instance, each sensor can update its estimate and then share the latter with its neighbors.

The device at each node of the sensor network is typically low-cost and battery-supplied. The latter feature, indeed, guarantees a certain ability of the device to adapt to the surrounding environment. On the other hand, the data  transmission represents the most energy consuming node task. Thus, it is fundamental to keep under control  data transmission which can be done by means of data-driven (or event-triggered) strategies for scheduling data communication, see \cite{BATTISTELLI2012926,7047754,shi2016event},
\cite{Liu20152470,li2015event,yan2014distributed}. In this paper we focus our attention to the distributed scheme proposed in \cite{BATTISTELLI201875}: each node updates its estimate with the new measurement (if available); then, the latter is compared with the one just propagated in time through the state space model (i.e. the estimate that can be computed also by the neighbors of that node in the case there is no transmission). If the discrepancy between them is large, then the node sends out the updated estimate to its neighbors. Finally, each node performs a fusion of its estimate and the ones from its neighbors. 
The appealing property of this transmission rule is that it allows to obtain a distributed algorithm which enjoys nice stability properties (i.e., mean-square boundedness of the state estimation error in all nodes) under minimal requirements.

In many situations the actual model is different from the nominal one (i.e. the one used in the estimation algorithm), \cite{RS_MPC_IET}. As a consequence, the resulting performance will be poor. One possible way to address this issue is to consider the robust Kalman filter proposed in  \cite{ROBUST_STATE_SPACE_LEVY_NIKOUKHAH_2013}: the idea is to consider a dynamic minimax game where one player is the estimator minimizing the state prediction error, while the other one is the hostile player which selects the least favorable model in a set of plausible models called ambiguity set. The latter is a ball which is formed by placing an upper bound on the Kullback-Leibler divergence between the nominal state space model and the models inside in. The radius of this ball is called tolerance and defines the magnitude of the uncertainty in the nominal model.

Distributed strategies which are robust to model uncertainty  have been already proposed in the literature, see for instance \cite{SHEN20101682,Luo_2008,RKDISTR,s20113244}, however, to the best of the authors' knowledge, none of them consider the case with event-triggered communication.

The contribution of this paper is to propose a new distributed Kalman filter with event-triggered communication under model uncertainty. More precisely, our approach represents a robust version of the distributed strategy proposed by \cite{BATTISTELLI201875}. Moreover, we also present a preliminary numerical example in order to test the performance of the proposed strategy.

The outline of the paper is as follows. In Section \ref{sec:formulation} we formulate the distributed state estimation problem characterized by data transmission constraints and model uncertainty. In Section \ref{sec:algo} we introduce the proposed robust distributed approach with event-triggered communication. In Section \ref{sec:numerical} we consider a numerical example showing the strength of our method. Finally, in Section \ref{sec:conclusion} we draw the conclusions.

 \section{Problem formulation}\label{sec:formulation}
Consider a network of nodes described by the digraph $(\mathcal{N}, \mathcal{A}, \mathcal{S})$ where: $\mathcal{N} = \{1,...,N\}$ is the set of the nodes, $\mathcal{A} \subseteq \mathcal{N}\times\mathcal{N}$ is the set of edges and $\mathcal{S} \subseteq \mathcal{N}$ is the subset of sensors node. The nodes in $\mathcal S$ are the only one that have the capabilities to perform  measurements, while the simple nodes in $\mathcal N\setminus \mathcal S$  are used to increase the connectivity of the network. If $(j,i)\in\mathcal A$, then it means that node $j$ can transmit data to node $i$; moreover, all the possible self-loops belong to $\mathcal A$. For each node $i\in\mathcal{N}$, the subset $\mathcal{N}_{i} := \{j:(j,i)\in\mathcal{A},\; j\neq i  \}\subseteq\mathcal{N}$ denotes the set of its in-neighbors, i.e. the set composed by the nodes that send information to node $i$.  

We attach to this network the following nominal state space model
\begin{align}
      x_{t+1} &= Ax_{t}+B\omega_{t}\label{eq:nominalStateEq}\\
    y^{i}_{t} & = C^{i}x_{t}+D^iv^{i}_{t}\;,\;\;\;\;  \;i\in\mathcal{S} \label{eq:nominalMeasEq}
\end{align}
where $x_{t}\in\mathbb R^n$ is the state, $y^i_t\in\mathbb R^{p_i}$ is the output at the sensor node $i$. Furthermore, matrix $B$ and $D^i$ are full row rank matrices, $\omega_{t}$ and $v^{i}_{t}$ are zero-mean normalized Gaussian white noises. The initial state $x_0$ is with mean $\mu_0$ and with covariance matrix
$P_0$. Finally, we assume that $\omega_{k}$'s, $v^{i}_{k}$'s and $x_0$ are independent. Notice that, the model (\ref{eq:nominalStateEq})-(\ref{eq:nominalMeasEq}) can be written as 
\begin{align}
      x_{t+1} &= Ax_{t}+B\omega_{t}\nn\\
    y_{t} & = Cx_{t}+D^iv_{t} \nn
\end{align}
where $C:=\col{C^i, i\in \mathcal S}$, i.e. it is the matrix obtained by stacking $C^i$'s, $D:=\mathrm{diag}(D^i, i\in \mathcal S)$ ,  i.e. it is the block diagonal matrix whose main blocks are $D^i$'s, $y_t:=\col{y_t^i, i\in \mathcal S}\in\mathbb R^p$ and $v_t:=\col{v_t^i, i\in \mathcal S}\in\mathbb R^p$.

In this paper we face a distributed state estimation problem under model uncertainty and data transmission constraints. More precisely, each node $i\in\mathcal N$ must estimate the state $x_t$ at each time instant $t\in\mathbb{Z}_{+}=\{1, ..., T\}$ given $Y_{t}:=\{y_0\ldots y_{t}\}$ and taking into account that: 
i) the actual model does not coincide with the nominal one (\ref{eq:nominalStateEq})-(\ref{eq:nominalMeasEq}); ii) each node $i$ can selectively transmits only the most relevant data, without compromising stability properties.

Regarding model uncertainty, we consider the framework proposed in \cite{ROBUST_STATE_SPACE_LEVY_NIKOUKHAH_2013,STATETAU_2017,ZORZI2017133} which is briefly reviewed below. The nominal model (\ref{eq:nominalStateEq})-(\ref{eq:nominalMeasEq}) over the time interval $t=1\ldots N$ can be equivalently described by the conditional probability densities $\phi_t(z_t|x_t)$, where $z_t=[\,x_{t+1}\; y_t\,]^T$, $t=1\ldots N$, and $f(x_0)$. We assume that the actual model is described by the conditional probability densities $\tilde \phi_t(z_t|x_t)$, $t=1\ldots N$, and $f(x_0)$. Moreover, we assume that $\tilde \phi_t$ belongs to the ambiguity set 
$\mathcal B_t=\{\,\tilde \phi_t \hbox{ s.t. } \tilde{\mathbb E}[\log(\tilde\phi_t/\phi_t) |Y_{t-1}]\leq b\,\}$,
\begin{align} &\tilde{\mathbb E}[\log(\tilde\phi_t/\phi_t) |Y_{t-1}]:=\nn\\
&\int_{\mathbb R^n}\int_{\mathbb R^{n+p}} \tilde \phi_t(z_t|x_t) \tilde f_t(x_t|Y_{t-1})\log\left( \frac{\tilde \phi_t(z_t|x_t) }{\phi_t(z_t|x_t) }\right)\mathrm d z_t \mathrm d x_t,\end{align}
$\tilde f_t$ is the actual probability density of $x_t$ given $Y_{t-1}$ and $b>0$ is called tolerance which accounts for model uncertainty. In plain words, $\mathcal B_t$ is a ball about the nominal density with with radius $b$.

\cite{ROBUST_STATE_SPACE_LEVY_NIKOUKHAH_2013} proposed as robust state estimator the one solving the following minimax game:
\begin{align} x_{t+1|t}=\underset{g_t\in\mathcal G}{\argmin}\underset{\tilde \phi_t\in\mathcal B_t}{\max}\, \tilde{\mathbb E}[\| x_{t+1}-g_t(y_t)\|^2|Y_{t-1}]\label{central_pb}
\end{align}
where  $ x_{t+1|t}$ is the estimator of $x_{t+1}$ given $Y_{t}$; $\mathcal G$ is the set of estimators having finite second order moments for any $\tilde \phi_t\in \mathcal B_t$; \begin{align}
    \tilde{\mathbb{E}}&\left [ \|x_{t+1}-g_{t}(y_{t})\|^{2}|Y_{t-1} \right] =\nn\\ & \int_{\mathbb R^{n}}\int_{\mathbb R^{n+p}}\|x_{t+1}-g_{t}(y_{t})\|^{2}\tilde{\phi}_{t}(z_{t}|x_{t})\tilde{f}_{t}(x_{t}|Y_{t-1})dz_{t}dx_{t}
\end{align}
and it is assumed that $\tilde{f}_{t}(x_{t}|Y_{t-1})\sim\mathcal{N}(\hat{x}_{t|t-1}, V_{t|t-1})$. The basic idea behind this paradigm is that when we are looking for an estimator that minimize properly the selected loss function, a hostile player called ``nature'' conspires to select the worst possible model in the ambiguity set $\mathcal{B}_{t}$. In \cite{ROBUST_STATE_SPACE_LEVY_NIKOUKHAH_2013} it has been shown that the (centralized) robust estimator solution to (\ref{central_pb}) admits a Kalman-like structure. It is not difficult to show that such a filter can be written in the information form as follows. Let $P_{t|t-1}$ and $V_{t|t-1}$ denote the pseudo-nominal and least favorable, respectively, covariance matrix of the prediction error at time $t$; let $x_{t|t}$ denote the estimator of $x_t$ given $Y_t$ and $P_{t|t}$ denotes the covariance matrices of the corresponding estimation error at time $t$.  We define the corresponding matrices in the information form as $\Omega_{t|t-1} = P_{t|t-1}^{-1}$, $\Psi_{t|t-1} = V_{t|t-1}^{-1}$, $\Omega_{t|t} = P_{t|t}^{-1}$ and the information states as
\begin{equation*}
q_{t|t-1}=\Psi_{t|t-1}{x}_{t|t-1},\quad q_{t|t}=\Omega_{t|t}{x}_{t|t}.
\end{equation*} Then, it is not difficult to see that the robust estimator obeys to
\begin{align}
      {\small \text{Correction step:}}&
    \begin{cases}\label{eq:RKFfiltering}
        \Omega_{t|t} = \Psi_{t|t-1}+ C^{T}R^{-1}C\\
        q_{t|t} = q_{t|t-1} + C^{T}R^{-1}y_{t} ,      
    \end{cases}\\
       {\small\text{Prediction step:}}&
    \begin{cases}\label{eq:RKFprediction}
    \Omega_{t+1|t} =\\ \hspace{0.2cm}     Q^{-1}-Q^{-1}A(A^{T}Q^{-1}A+\Omega_{t|t})A^{T}Q^{-1}\\
                \text{Find } \theta_{t}>0 \text{ s.t. }\gamma(\Omega_{t+1|t}, \theta_{t})=b\\
  \Psi_{t+1|t} = \Omega_{t+1|t}-\theta_{t}\mathbf{I}_{n}\\
        q_{t+1|t} = \Psi_{t+1|t}A\Omega^{-1}_{t|t}q_{t|t},    
    \end{cases}
\end{align}
 where $R:=DD^{T}$, $Q:=BB^{T}$ are positive definite matrices and
 \al{ \gamma&(\Omega,\theta) \nn\\ &:= \frac{1}{2}\left\{ tr[(\mathbf{I}_{n}-\theta \Omega^{-1})^{-1}-\mathbf{I}_{n}] + \log\det(\mathbf{I}_{n}-\theta \Omega^{-1})\right \}.\nn }The parameter $\theta_{t}>0$ is called risk sensitivity parameter. It is worth noting that given $\Omega>0$ and $b>0$, the equation $\gamma(\Omega,\theta)=b$ always admits a unique solution $\theta>0$ such that $\Omega-\theta \mathbf{I}_{n}>0$. Furthermore, in the special case where $b =0$, i.e. there is no model uncertainty, then $\gamma(\Omega,\theta)=0$ implies that $\theta_{t}=0$ and thus the above equations degenerates in the usual Kalman equations in the information form.

  \section{Proposed algorithm}\label{sec:algo}
 Before to introduce our robust distributed estimation paradigm we consider the following quite simple scenario. We assume that $\mathcal N=\mathcal S$, i.e. all the nodes are sensor nodes, and $\mathcal A=\emptyset$, i.e. the nodes do not communicate. In the presence of model uncertainty, at node $i$ we can consider the robust Kalman filter in the information form based on the local model
 \begin{align}
      x_{t+1} &= Ax_{t}+B\omega_{t}\label{eq:nominalStateEqloc}\\
    y^{i}_{t} & = C^{i}x_{t}+D^iv^{i}_{t}\; ;   \label{eq:nominalMeasEqloc}
\end{align}
thus, we obtain the following algorithm:
\begin{align}
     {\small \text{Correction step:}}&
    \begin{cases}\label{eq:RKFfilteringloc}
        \Omega_{t|t}^i = \Psi_{t|t-1}^i+ (C^i)^{T}(R^i)^{-1}C^i\\
        q_{t|t}^i = q_{t|t-1}^i + (C^i)^{T}(R^i)^{-1}y_{t}^i       
    \end{cases}\\
    {\small \text{Prediction step:}}&
    \begin{cases}\label{eq:RKFpredictionloc}
    \Omega_{t+1|t}^i =\\ \hspace{0.2cm}     Q^{-1}-Q^{-1}A(A^{T}Q^{-1}A+\Omega_{t|t}^i)A^{T}Q^{-1}\\
                \text{Find } \theta_{t}^i>0 \text{ s.t. }\gamma(\Omega_{t+1|t}^i, \theta_{t}^i)=b\\
  \Psi_{t+1|t}^i = \Omega_{t+1|t}^i-\theta_{t}^i\mathbf{I}_{n}\\
        q_{t+1|t}^i= \Psi_{t+1|t}^iA(\Omega^{i}_{t|t})^{-1}q_{t|t}^i
    \end{cases}
\end{align}
 where $R^i:=D^i(D^i)^{T}$. It is worth noting that each node has its own risk sensitivity parameter $\theta_{t}^i$. 
  
Next, we consider the scenario in which: i) $ \mathcal S\subseteq \mathcal N$; ii) each node $i$ can transmit its local estimate ${q}^{i}_{t|t}$ and information matrix $\Omega^{i}_{t|t}$ to all its out-neighbors $\{\,j\, : \, (i,j)\in\mathcal A, \, i\neq j\,\}$ if necessary.  In plain words, each node $i$ can decide at any time step whether to transmit or not its data, i.e. $(q^i_{t|t},\Omega^{i}_{t|t})$,  without compromising the stability properties of the algorithm, i.e. $\|x_t-x^i_{t|t}\|$ and $\|P^{i}_{t|t}\|$ do not diverge for any $i\in\mathcal N$ as $t$ approaches infinity.

The estimation paradigm that we now present  
is a robust extension of the distributed state estimation algorithm with event-triggered communication proposed in  \cite{BATTISTELLI201875} and it is composed by four steps described below.

\textbf{Correction.} At time $t$,  the predicted pair $(q^{i}_{t|t-1}, \Psi^{i}_{t|t-1})$ is available at node $i\in \mathcal N$. If $i\in\mathcal S$, i.e. it is a sensor node, then also the measurement $y_t^i$ is available and thus the correction step coincides with (\ref{eq:RKFfilteringloc}). If $i\notin\mathcal S$,
no measurement is available at the node, then we can only   propagate the prediction couple. Therefore, the so called information pair is obtained as
\begin{equation}\label{eq:correctionStep}
    (q^{i}_{t|t},\Omega^{i}_{t|t}) =
    \begin{cases}
        \text{use equations (\ref{eq:RKFfilteringloc}),} & \text{if $i\in\mathcal{S}$}\\
        (q^{i}_{t|t-1}, \Psi^{i}_{t|t-1}), & \text{if $i\in\mathcal{N \setminus 
        S}$}.
    \end{cases}
\end{equation}

\textbf{Information exchange.} Each node $i\in\mathcal N$ sends its information couple $(q^{i}_{t|t}, \Omega^{i}_{t|t})$ to all its out-neighbors according to the binary variable  $c^{i}_{t}$: \begin{itemize} 
\item if $c^{i}_{t}=1$, then node $i $ transmits the information couple to all its out-neighbors at time $t$;
 \item  if $c^{i}_{t}=0$, then node $i $ does not transmit the  information couple to all its out-neighbors at time $t$.
 \end{itemize}
 It remains to define the binary variable $c^i_t$. Let $n^{i}_{t}\in\mathbb N$ be the number of time instants elapsed from the most recent transmission of node $i$, i.e.  the most recently transmitted data is $(q^{i}_{t-n^{i}_{r}|r-n^{i}_{t}}, \Omega^{i}_{t-n^{i}_{t}|t-n^{i}_{t}})$. Then, all the out-neighbors of node $i$ propagate $(q^{i}_{t-n^{i}_{r}|t-n^{i}_{t}}, \Omega^{i}_{t-n^{i}_{t}|t-n^{i}_{t}})$ in time through a prediction step which takes into account the fact that the actual model does not coincide with the nominal one (see (\ref{eq:estPredictionStep}) in the prediction step below). Let  $(\bar{q}^{i}_{t}, \bar{\Omega}^{i}_{t})$ denote this propagated pair at time $t$. Then, the transmission rule $c_t^i$ computed at node $i$ is defined as in \cite{BATTISTELLI201875}:  
 \begin{equation}\label{eq:eventTriggered}
    c^{i}_{t} = 
    \begin{cases}
        0, & \quad \text{if } \| {x}^{i}_{t|t}-\bar{x}^{i}_{t}\|^{2}_{\Omega^{i}_{t|t}}\leq\alpha\\ & \text{ and } \frac{1}{1+\beta}\Omega^{i}_{t|t}\leq \bar{\Omega}^{i}_{t} \leq (1+\delta)\Omega^{i}_{t|t}\\
        1, & \quad \text{otherwise}
    \end{cases}
\end{equation}
where $\bar{x}^{i}_{t} = (\bar{\Omega}^{i}_{t})^{-1}\bar{q}^{i}_{t}$ represents the state prediction based on the propagation of the most recent transmitted pair  $(q^{i}_{t-n^{i}_{r}|t-n^{i}_{t}}, \Omega^{i}_{t-n^{i}_{t}|t-n^{i}_{t}})$; $x^{i}_{t|t} = ({\Omega}^{i}_{t|t})^{-1}q^{i}_{t|t}$ is the state estimate at node $i$. In plain words, the transmission rule in (\ref{eq:eventTriggered})
checks the discrepancy between $(q^{i}_{t|t}, \Omega^{i}_{t|t})$ and $(\bar{q}^{i}_{t}, \bar{\Omega}^{i}_{t})$. If the latter is large,
then it means that the out-neighbors own a prediction corresponding to node $i$ which is bad and thus node $i$ must transmit the data.   
The positive scalars $\alpha$, $\beta$ and $\delta$ can be tuned by the user in order to reach a desired behavior in terms of transmission rate and performance. More precisely, $\alpha$ tunes the bound on the discrepancy between $x^i_{t|t}$ and $\bar x^i_t$, while $\beta$ and $\delta$ tunes the allowed mismatch between the covariance matrices ${\Omega}^i_{t|t}$ and $\bar{\Omega}^i_t$. In   \cite{BATTISTELLI16} it has been shown that   
the transmission strategy in (\ref{eq:eventTriggered}) guarantees the following upper bound. If we model the propagated and the information pairs as $\mathcal N(\bar{q}^{i}_{t}, \bar{\Omega}^{i}_{t})$ and $\mathcal N (q^{i}_{t|t}, \Omega^{i}_{t|t})$,  respectively,  then  \begin{equation*} 
    D_{KL}(\mathcal N(q^{i}_{t|t}, \Omega^{i}_{t|t}),\mathcal N(\bar{q}^{i}_{t}, \bar{\Omega}^{i}_{t})) \leq \frac{1}{2}[\alpha+\beta n+n\log(1+\delta)]
\end{equation*}
where $n$ is the state dimension and $D_{KL}$ denotes the Kullback-Leibler divergence.

\textbf{Information fusion.} In this step, any node merges its information with the ones regarding its in-neighbors. Let $\Pi\in\mathbb R^{N\times N}$ denote the consensus matrix whose element
in position $(i,j)$ is defined as:
\begin{equation*}
    \pi_{i,j} = \begin{cases}
          (d_{i}+1) ^{-1}, & \text{ if }(j,i)\in\mathcal{A}\\
          0, & \text{otherwise}
    \end{cases}
\end{equation*} where $d_i$ denotes the degree of node $i$; in this way we have that 
$\pi_{i,j}$ with $j\in\mathcal N_i$ represents the coefficients of a convex combination. 
 Then, the fusion step is performed through the following convex combination of the  pairs:
 \begin{equation}\label{eq:qFusedwithTilde}
    q^{i,F}_{t|t} = \pi_{i,i}q^{i}_{t|t}+\sum_{j\in\mathcal{N}_{i}}\pi_{i,j}\left [c^{j}_{t}q^{j}_{t|t}+(1-c^{j}_{t})\tilde{q}^{j}_{t}\right ]
\end{equation}
\begin{equation}\label{eq:OmegaFusedwithTilde}
    \Omega^{i,F}_{t|t} = \pi_{i,i}\Omega^{i}_{t|t}+\sum_{j\in\mathcal{N}_{i}}\pi_{i,j}\left [c^{j}_{t}\Omega^{j}_{t|t}+(1-c^{j}_{t})\tilde{\Omega}^{j}_{t}\right ] 
\end{equation}
where 
\begin{equation*}
    \tilde{q}^{j}_{t} = \frac{1}{1+\delta}\bar{q}^{j}_{t}\;\;\;,\;\;\;\tilde{\Omega}^{j}_{t} = \frac{1}{1+\delta}\bar{\Omega}^{j}_{t}.
\end{equation*}
In view of (\ref{eq:qFusedwithTilde})-(\ref{eq:OmegaFusedwithTilde}), we can see that in the fusion step we consider $(q^j_{t|t},\Omega^j_{t|t})$, if node $j$ transmitted its information pair at time $t$. If node $j$ does not transmit, then the aforementioned pair is not available at node $i$. To account for this lack, given that at each iteration the nodes can calculate the pair $(\bar{q}^{i}_{t}, \bar{\Omega}^{i}_{t})$, which is certainly less informative than $(q^{i}_{t|t}, \Omega^{i}_{t|t})$, it becomes convenient to shrink  $(\bar{q}^{i}_{t}, \bar{\Omega}^{i}_{t})$ in  (\ref{eq:qFusedwithTilde})-(\ref{eq:OmegaFusedwithTilde}) by   the factor $(1+\delta)^{-1}$ in order to decrease its importance in the fusion step, see \cite{BATTISTELLI201875} for more details.

\textbf{Prediction.} Once each node $i\in\mathcal{N}$ has computed the fused information couple $(q^{i,F}_{t|t} \Omega^{i,F}_{t|t})$, the latter is propagated in  time with the robust prediction step in (\ref{eq:RKFpredictionloc}) where $(q^i_{t|t},\Omega_{t|t}^{i})$ is now replaced by $(q^{i,F}_{t|t},\Omega_{t|t}^{i,F})$:
\begin{align}\label{eq:infoFormPred} 
\left.\begin{array}{ll}
    \Omega_{t+1|t}^i &=Q^{-1}-Q^{-1}A(A^{T}Q^{-1}A+\Omega_{t|t}^{i,F})A^{T}Q^{-1}\\
                \text{Find } &\theta_{t}^i>0 \text{ s.t. }\gamma(\Omega_{t+1|t}^i, \theta_{t}^i)=b\\
  \Psi_{t+1|t}^i &= \Omega_{t+1|t}^i-\theta_{t}^i\mathbf{I}_{n}\\
        q_{t+1|t}^i&= \Psi_{t+1|t}^iA(\Omega^{i,F}_{t|t})^{-1}q_{t|t}^{i,F}.
\end{array}\right. 
\end{align}
In this step we also need to propagate in time
the pair $(\bar{q}^{i}_{t}, \bar{\Omega}^{i}_{t})$, i.e. the one used in the case node $i$ does not transmit its information pair. 
Notice this operation is performed by both node $i$ and its out-neighbors and it can be summarized as follows.  At nodes $i\,\cup\, \mathcal N_i $  we have the pair 
$(\breve{q}^{i}_{t},\breve{\Omega}^{i}_{t})$ defined as
\begin{align}
\begin{cases}\label{def_breve}
\breve{q}^{i}_{t} &= c^{i}_{t}q^{i}_{t|t}+(1-c^{i}_{t})\bar{q}^{i}_{t}\\ 
\breve{\Omega}^{i}_{t} &= c^{i}_{t}\Omega^{i}_{t|t}+(1-c^{i}_{t})\bar{\Psi}^{i}_{t};
\end{cases}
\end{align} then, it is propagated in  time with  the robust prediction step in (\ref{eq:RKFpredictionloc}) where $(q^i_{t|t},\Omega_{t|t}^{i})$ is now replaced by $(\breve q^{i}_{t},\breve \Omega_{t}^{i})$: 
\begin{align}
\label{eq:estPredictionStep}
\left.\begin{array}{ll}
   \bar{ \Omega}_{t+1}^i &=Q^{-1}-Q^{-1}A(A^{T}Q^{-1}A+\breve\Omega_{t}^{i})A^{T}Q^{-1}\\
                \text{Find } &\bar\theta_{t}^i>0 \text{ s.t. }\gamma(\bar \Omega_{t+1}^i, \bar\theta_{t}^i)=b\\
 \bar \Psi_{t+1}^i &= \bar\Omega_{t+1}^i-\bar\theta_{t}^i\mathbf{I}_{n}\\
        \bar q_{t+1}^i&= \bar \Psi_{t+1}^iA(\breve\Omega^{i}_{t})^{-1}\breve q_{t}^{i}.
\end{array}\right. 
\end{align}

The procedure is summarized in Algorithm \ref{alg:DKFeventTriggered} below. It is interesting to note that each node $i$ is characterized by two risk sensitivity parameters, i.e. $\theta^i_t$ and $\bar \theta^i_t$. In the case that $b=0$, i.e. there is no model uncertainty, in Algorithm \ref{alg:DKFeventTriggered} we have: $\theta^i_t=0$, $\bar \theta^i_t=0$ and thus $\Psi_{t}^i = \Omega_{t}^i$, $ \bar \Psi_{t}^i = \bar\Omega_{t}^i$, i.e. we recover the distributed Kalman algorithm with event-triggered communication proposed in \cite{BATTISTELLI201875}.
\begin{algorithm}
    \caption{\textbf{RDKF} with event-triggered communication}\label{alg:DKFeventTriggered}
    \textbf{Initialization:} Set $(q^{i}_{0|-1}, \Psi^{i}_{0|-1})$ for any $i\in\mathcal N$
    \newline
    For each $t=0,1,\ldots$
    \newline
    For each node $i\in\mathcal{N}$
    \begin{description}
        \item \textbf{Correction:} 
        \begin{equation*}
    (q^{i}_{t|t},\Omega^{i}_{t|t}) =
    \begin{cases}
        \text{use   (\ref{eq:RKFfilteringloc}),} & \text{if $i\in\mathcal{S}$}\\
        (q^{i}_{t|t-1}, \Psi^{i}_{t|t-1}), & \text{if $i\in\mathcal{N \setminus 
        S}$}
    \end{cases}
\end{equation*}
        \item \textbf{Information exchange:} 
        \begin{description}
        \item[-] if $t=0$ set $c^i_t=1$, otherwise determine $c_t^i$ according to (\ref{eq:eventTriggered})
        \item[-] if $c^i_t=1$ transmit $(q^i_{t|t},\Omega^i_{t|t})$ to the out-neighbors 
        \item[-] receive $(q^j_{t|t},\Omega^j_{t|t})$ from all the in-neighbors $j\in\mathcal N_i$ for which $c_t^j=1$
        \end{description}
\item \textbf{Information fusion:} 
\begin{equation*}
\hspace*{-0cm}    \tilde{q}^{j}_{t} = \frac{1}{1+\delta}\bar{q}^{j}_{t}\;\;\;,\;\;\;\tilde{\Omega}^{j}_{t} = \frac{1}{1+\delta}\bar{\Omega}^{j}_{t} \text{ with } j\in\mathcal N_i
\end{equation*}
 \begin{align*}
    q^{i,F}_{t|t} &= \pi_{i,i}q^{i}_{t|t}+\sum_{j\in\mathcal{N}_{i}}\pi_{i,j}\left [c^{j}_{t}q^{j}_{t|t}+(1-c^{j}_{t})\tilde{q}^{j}_{t}\right ]\\
    \Omega^{i,F}_{t|t} &= \pi_{i,i}\Omega^{i}_{t|t}+\sum_{j\in\mathcal{N}_{i}}\pi_{i,j}\left [c^{j}_{t}\Omega^{j}_{t|t}+(1-c^{j}_{t})\tilde{\Omega}^{j}_{t}\right ] 
\end{align*}
\item \textbf{Prediction step:} 
 \begin{description}
 \item[-] Compute $(q^i_{t+1|t},\Omega^i_{t+1|t})$ using (\ref{eq:infoFormPred})
 \item[-] Compute $(\bar q^i_{t+1|t},\bar \Psi^i_{t+1|t})$ using (\ref{def_breve}) and (\ref{eq:estPredictionStep})
 \end{description}
     \end{description}
\end{algorithm}


 \section{Numerical example}\label{sec:numerical}
  \begin{figure}
\begin{center}
\includegraphics[width=0.5\textwidth]{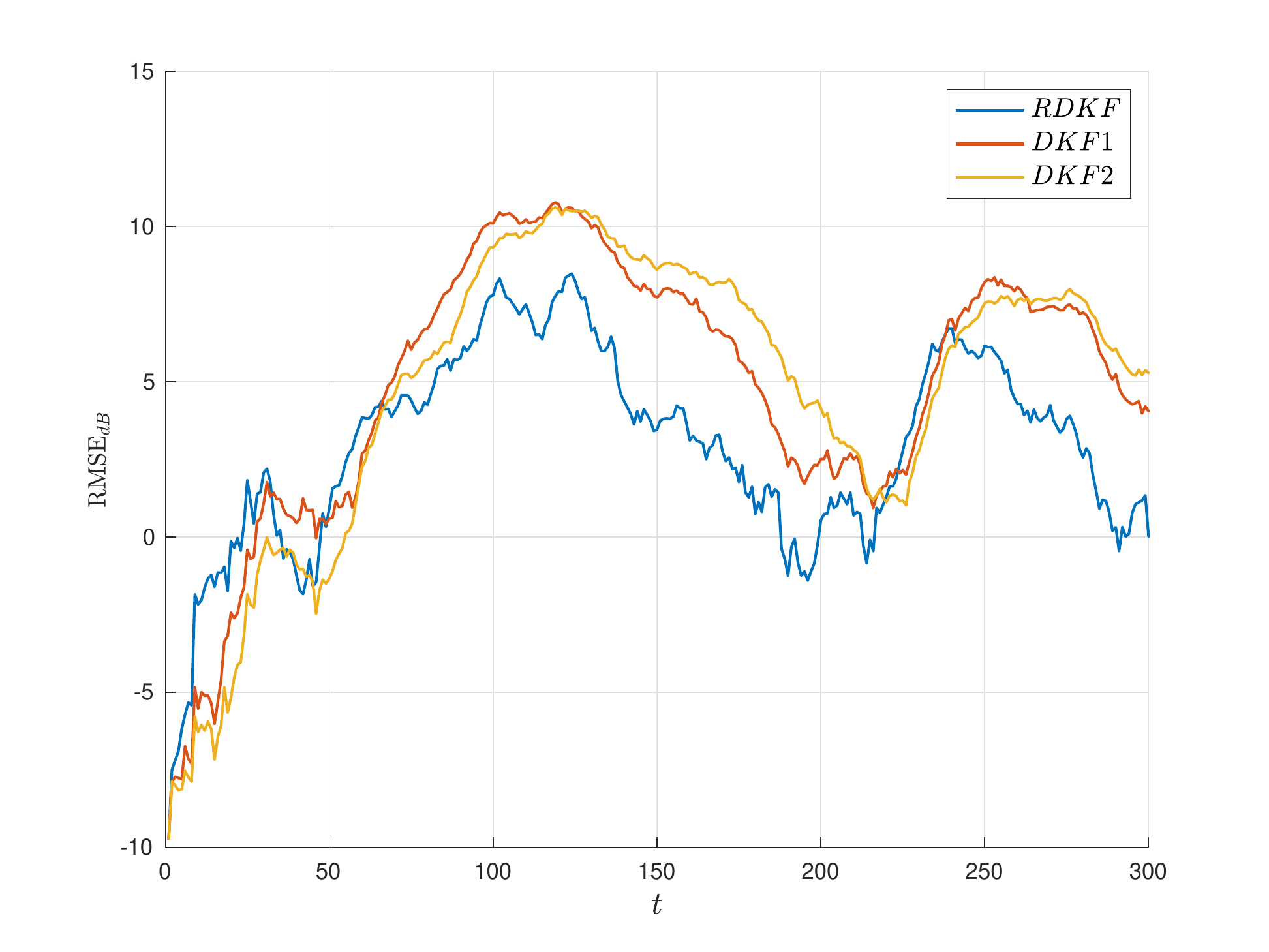}    
\caption{Average root mean square error across the network.} 
\label{fig:rmse}
\end{center}
\end{figure}

 \begin{figure*}
\begin{center}
\includegraphics[width=\textwidth]{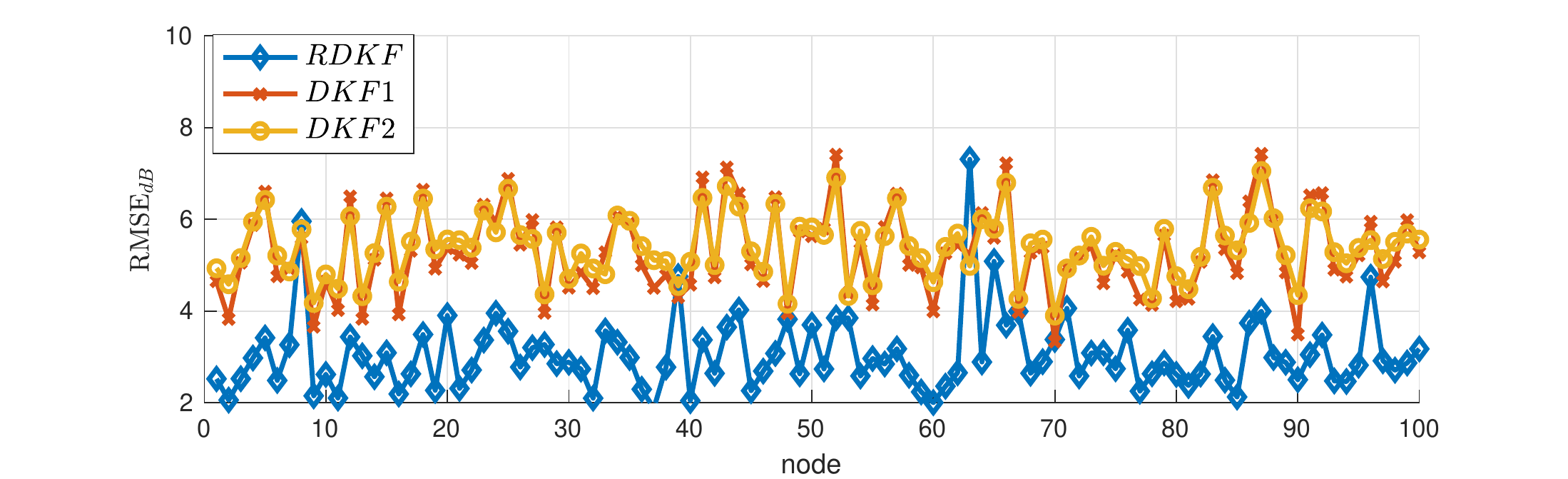}    
\caption{Average root mean square error at each node over the time horizon $[1,300]$.} 
\label{fig:rmsen}
\end{center}
\end{figure*}

 \begin{figure}
\begin{center}
\includegraphics[width=0.5\textwidth]{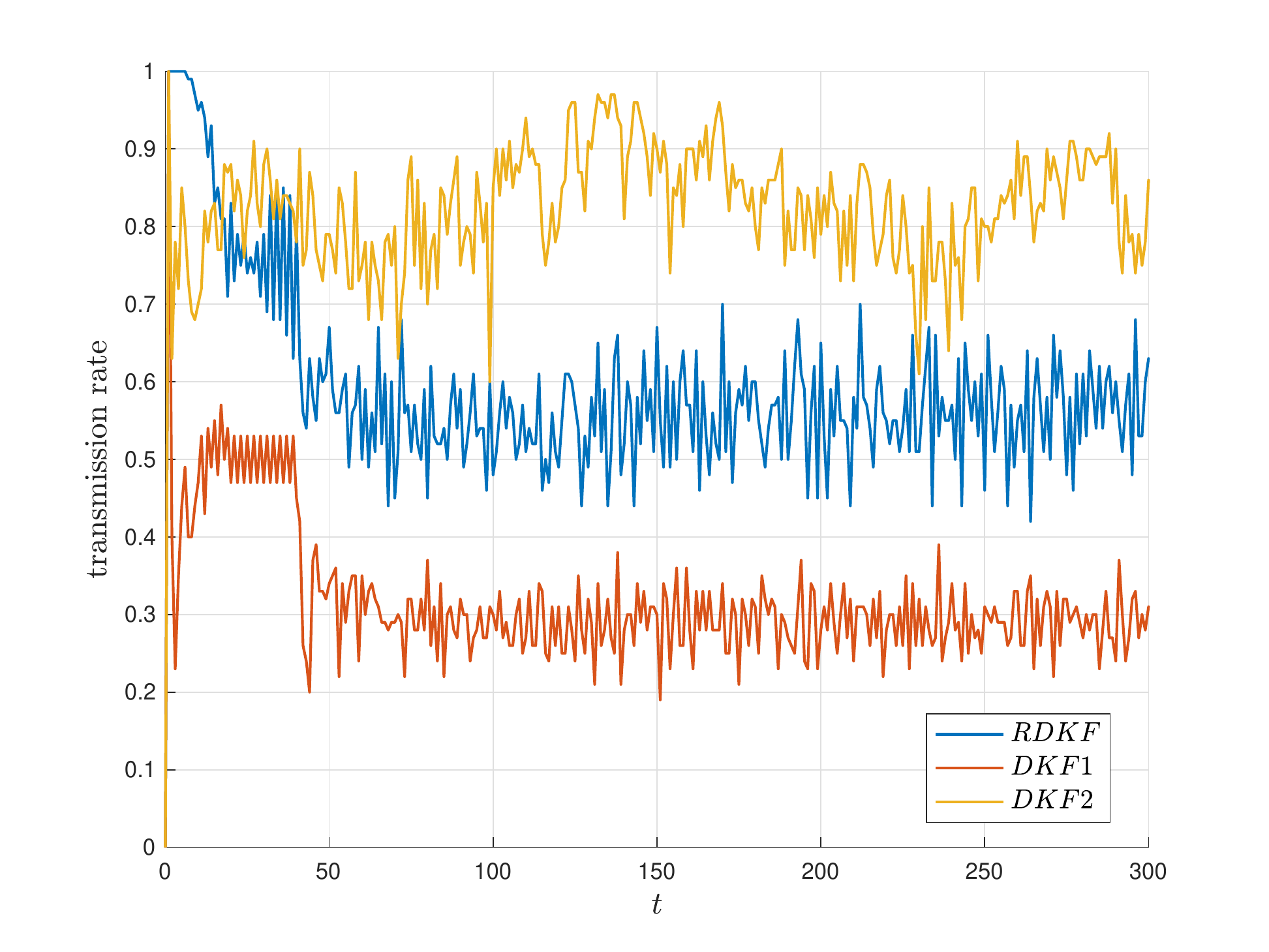}    
\caption{Transmission rate across the network.} 
\label{fig:erre}
\end{center}
\end{figure}

 \begin{figure}
\begin{center}
\includegraphics[width=0.5\textwidth]{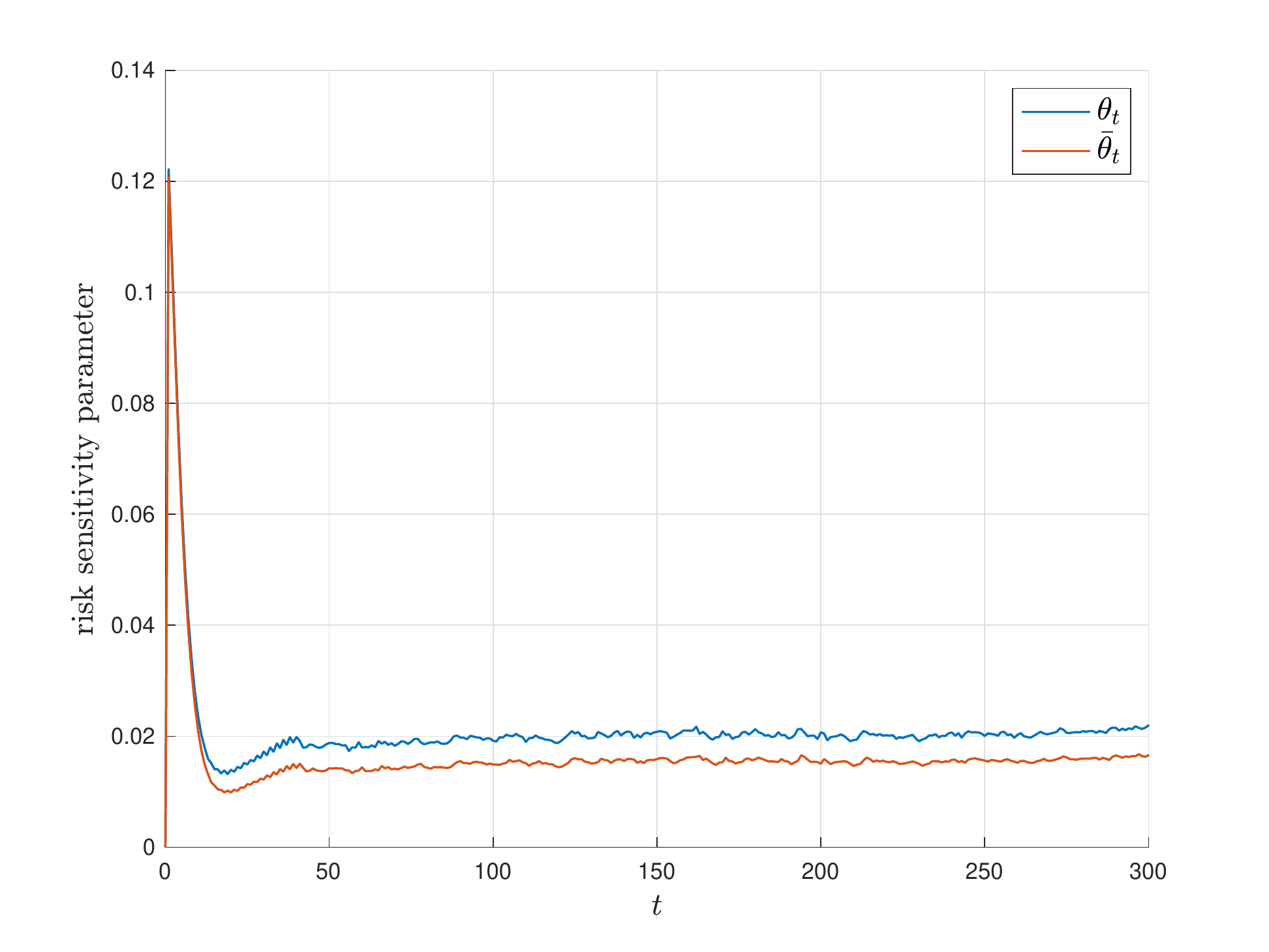}    
\caption{Average risk sensitivity parameters across the network.} 
\label{fig:theta}
\end{center}
\end{figure}
In this section we evaluate the performance of the proposed robust distributed Kalman algorithm with event-triggered communication. We consider the problem in \cite{cattivelli2010} of tracking the position of a projectile by using noisy position measurements obtained by a network of $N=100$ nodes where $20$ of them are sensor nodes. The possible connections among the nodes has been randomly generated in such a way the network is strongly connected. The model for the projectile motion is 
\al{\label{proj_cont}\dot x^c_t =\Phi x^c_t+u^c_t}
where 
\al{\Phi=\left[\begin{array}{cc}0 & 0 \\ \mathbf I_3 & 0\end{array}\right],\nn}
$u^c_t=[\, 0 \;0 \; -g\; 0 \; 0\; 0\,]^T$, with $g=-10$, and $x^c_t=[\, v_{x,t}\;v_{y,t}\;v_{z,t}\;p_{x,t}\;p_{y,t}\;p_{z,t}\;\,]^T$ with $v$ denoting the velocity, $p$ the position and the subscripts $x,y,z$ denoting the three spatial dimensions. We discretize (\ref{proj_cont}) with sampling time equal to $0.1$. The corresponding discrete time model is $x_{t+1}=Ax_{t}+u_t$ where $x_t$ is the sampled version of $x^c_t$, $A=\mathbf I_6+0.1\Phi$ and $u_t=(0.1 \mathbf I_6+0.1^2\Phi/2)u^c_t$. We assume that every sensor measures the position of the projectile in either two horizontal dimensions, or a combination of one horizontal dimension and the vertical dimension; in plain words, one sensor does not have measurements in all the three dimensions. Therefore, we obtain the  nominal discrete state-space model (\ref{eq:nominalStateEq})-(\ref{eq:nominalMeasEq}) where $C^i=[\, 0\; 0 \;0 \;\mathrm{diag(1,1,0)}\,]$, in the case that the sensor measures only the horizontal positions, or $C^i=[\, 0\; 0 \;0 \;\mathrm{diag(1,0,1)}\,]$, $C^i=[\, 0\; 0 \;0 \;\mathrm{diag(0,1,1)}\,]$, in the case that the sensor measures one horizontal position and the vertical position. Moreover, we choose $B=\sqrt{0.001}I$, $R^i=D^i(D^i)^T=\sqrt{k} PR_0P^T$ where $R_0=0.5\cdot \mathrm{diag}(1,4,7)$ and $P$ is a permutation matrix randomly chosen for every node. Finally, the initial state $x_0$ is a Gaussian random vector with covariance matrix $P_0=I$. Since the previous model is just an idealization of the underlying physical system, we assume that the actual state-space model belongs to the ambiguity set $\mathcal B_t$ about the aforementioned nominal model and with tolerance 
$b=0.05$. More precisely, we assume that the actual model is  the least favorable model solution to (i.e. the maximizer of) the centralized problem (\ref{central_pb}).

In the following we consider the following distributed filters:
\begin{itemize}
\item \textsf{RDKF} -- the distributed robust Kalman filter with event-triggered communication in Algorithm \ref{alg:DKFeventTriggered} with $b=0.05$. Here, the 
transmission rule (\ref{eq:eventTriggered}) is with $\alpha=30$, $\beta=0.2$ and $\delta=0.1$;
\item \textsf{DKF1} -- the distributed Kalman filter with event-triggered communication proposed in \cite{BATTISTELLI201875} and the transmission rule is with $\alpha=30$, $\beta=0.2$ and $\delta=0.1$;
\item \textsf{DKF2} -- is the same as \textsf{DKF1} but the transmission rule is with $\alpha=0.01$, $\beta=0.2$ and $\delta=0.1$.
\end{itemize}
It is worth noting that \textsf{RDKF} and \textsf{DKF1} have the same parameters for the transmission rule. However, taking the parameter $\alpha$ in (\ref{eq:eventTriggered}) the same for both \textsf{RDKF} and \textsf{DKF1} provides a transmission rate for \textsf{DKF1} which is smaller than the one of \textsf{RDKF} (see below). For this reason, we also consider \textsf{DKF2} where the parameter $\alpha$ has been decreased in oder to increase the transmission rate. We also tried to increase the transmission rate by keeping fixed $\alpha=0.01$ and changing $\beta,\delta$; however, we did not notice a significant growth in terms of transmission rate. 


We consider the following performance indexes:
\begin{itemize}
\item The average root mean square error across the network at time $t$: $$ \mathrm{RMSE}_t=\frac{1}{N}\sum_{i=1}^N \|x_t- x_{t|t}^i\|^2.$$
\item The average root mean square error at node $i$ over the time horizon $[1,300]$:
$$ \mathrm{RMSE}_i=\frac{1}{300}\sum_{t=1}^{300} \|x_t- x_{t|t}^i\|^2.$$ 
\item The transmission rate across the network at time $t$, i.e. the faction of nodes that transmit their data at time $t$. 
\end{itemize}
Figure \ref{fig:rmse} and  Figure \ref{fig:rmsen} show the two aforementioned root mean square errors, while Figure \ref{fig:erre} shows the transmission rate. As we can see  \textsf{RDKF}   outperforms \textsf{DKF1}. On the other hand,  the transmission rate of  \textsf{DKF1} is smaller than the one of \textsf{DKF}. However, even in the case we increase the transmission rate, i.e. we consider \textsf{DKF2}, \textsf{RDKF} is the best estimator.

 Finally, Figure \ref{fig:theta} shows the average risk sensitivity  parameters across the network:
 \al{\theta_t&=\frac{1}{N}\sum_{i=1}^N\theta_t^i,\quad 
 \bar \theta_t=\frac{1}{N}\sum_{i=1}^N	\bar \theta_t^i.\nn}
 We can notice that $\theta_t\geq \bar \theta_t$, moreover we have checked that $\theta_t^i\geq \bar \theta_t^i$ for many nodes. Since the mapping $\theta\mapsto \gamma( \Omega,\theta)$, with $\Omega>0$, is monotone increasing,  \cite{CONVTAU,ZORZI_CONTRACTION_CDC,LEVY_ZORZI_RISK_CONTRACTION}, and in view of the fact that
 $$   \gamma( \Omega_{t+1|t}^i,\theta_t^i)=\gamma( \bar \Omega^i_{t+1},\bar \theta_t^i)=b,$$
 it follows that $\Omega_{t+1|t}^i\geq  \bar \Omega^i_{t+1}$ for many nodes.


\section{Conclusion}\label{sec:conclusion}
We have proposed a robust distributed Kalman filter with event-triggered communication. More precisely, each sensor designs its predictor according to the least favorable model in a prescribed ambiguity set. The latter is ball about the nominal model. Moreover, each node transmits its data to its neighbors only in the case the discrepancy between the information pair, i.e. the one if we transmit, and the propagated pair, i.e. the one obtained without transmission, is not negligible. A numerical example showed the effectiveness of the proposed method in the case the actual model is different from the nominal one.

It remains to investigate the stability of the proposed distributed filter. We conjecture that, under the assumptions made in \cite{BATTISTELLI201875}, that is  
\begin{itemize}
  \item[\textbf{A1.}] the system matrix A is invertible,
  \item[\textbf{A2.}] the system is \textit{collectively observable}, i.e., the pair \textit{(A,C)} is observable where $C = \col{(C^{i};\;i\in\mathcal{S})}$,
  \item[\textbf{A3.}] the network is strongly connected, i.e. there exists a directed path between any pair of nodes $i,j \in \mathcal{N}$, 
\end{itemize}
the estimation error $x_t-x_{t|t}^i$ is bounded in mean-square  in each node.


\end{document}